\input amstex 
\documentstyle{amsppt}
\input bull-ppt
\keyedby{bull336e/mhm}

\topmatter
\cvol{28}
\cvolyear{1993}
\cmonth{January}
\cyear{1993}
\cvolno{1}
\cpgs{90-94}
\title A New Result for the Porous  Medium Equation  
\\Derived from the
Ricci Flow\endtitle
\author Lang-Fang Wu  \endauthor
\shortauthor{L.-F. Wu}
\shorttitle{A New Result for the Porous  Medium Equation}
\address Center for Math Analysis, 
Australian National  University, 
Canberra, Act 2601  \newline Australia\endaddress
\ml lang\@gauss.anu.edu.au\endml
\address {\it Current address}: Princeton University, 
Mathematics Department,
Princeton, NJ 08544\endaddress
\ml lfwu\@math.princeton.edu \endml
\date September 24, 1991 and, in revised form, on November 
13, 1991 and May 1,
1992. Presented on July 2 at the 1992 Regional Geometry 
Institute at Utah\enddate
\abstract Given $\Bbb R^2, $  with a ``good'' complete 
metric, we show that the
unique solution of the Ricci flow approaches a soliton at 
time infinity.
Solitons are solutions  of the Ricci flow, which move only 
by diffeomorphism.
The Ricci flow on $\Bbb R^2$ is the limiting case of the 
porous medium equation 
when $m$ is zero.  The results in the Ricci flow may 
therefore be interpreted
as sufficient conditions on the initial data, which 
guarantee that the
corresponding unique solution for the porous medium 
equation  on the entire
plane asymptotically behaves like a 
``soliton-solution''.\endabstract 
\endtopmatter

\document      

On $\Bbb R^2, $  any metric can be expressed as $ds^2 = 
e^u ( dx^2 + dy^2)$,
where $\{ x, y \}$ are rectangular coordinates on $\Bbb 
R^2$. Let $R$ be the
scalar curvature, then  the so-called Ricci flow on $\Bbb 
R^2$  is 
$$ \frac{\partial} {\partial t} ds^2 = -R ds^2, \tag 
"(*0)"$$
which may also be expressed as  
$$ \frac{\partial} {\partial t} u = e^{-u} 
\overline{\Delta} u, 
\quad \text {where }  
\overline{\Delta}=\partial_{x}^2 + \partial_{y}^2. \tag 
"(*1)"$$

The porous medium equation  is defined as 
$$\frac{\partial} {\partial \tau}  v = \overline{\Delta} 
v^m,\tag "(*2)" $$
where $ 0< m<\infty,$ and $v$ is a function on $\Bbb R^2$.
If we let $\tau= t/m$, then (*2) can be expressed as 
$$ \frac{\partial} {\partial t}  v = 
\frac{\overline{\Delta} (v^m -1)}{m}.$$
The limiting case of the porous medium equation as  $ m 
\to 0$ 
[PME $(m=0)]$ is therefore
 $$ \frac{\partial} {\partial t}  v =
\overline{\Delta} \ln v. \tag "(*3)"$$
In dimension 1, the PME $(m=0)$ with initial data $v 
\notin L^1$  has been
studied by  \cite {ERV,  H1,  H2, V}.

We are grateful to Sigurd B.~Angenent for pointing out the 
following observation
reflected in \cite {W2, Appendix}.

\thm{Proposition \cite {Angenent}} On $\Bbb R^2$, if we 
let  $v = e^u$, then
the limiting case of the porous medium  equation as $m \to 
0$ $(^*3)$ is
equivalent to the Ricci flow $(^*1)$.  \ethm

On a complete  $(\Bbb R^2, ds^2)$ we may define:

{\rm (a)} the {\it circumference at infinity\/} to be
$$C_{\infty}(ds^2) = {\underset K \to \sup} 
\, { \underset { D^2} \to \inf }\{  L( \partial
D^2)|   \forall \text{ compact set } K \subset \Bbb R^2,
\forall \text{ open set } D^2 \supset K\}$$ 
where $L(\partial D^2)$ is the length of $\partial D^2$ 
with respect to the
metric  $ds^2$;

(b) the {\it aperture\/} to be 
$$ A(ds^2)=
\frac{1}{2\pi}\lim_{r \to
\infty}\frac{L(\partial B_r)}{r},$$ 
 where       
 $B_r$ is a  geodesic ball at any given point on
$\Bbb R^2$   with
radius $r$.

A gradient  soliton  is a solution of the Ricci flow, 
which moves only
by 
diffeomorphism and there exists a function $f$ such that  $ 
\frac{\partial}{\partial t}g_{ij}  = L_{\bigtriangledown 
f} g_{ij},$ where
$L_{\bigtriangledown f}$ is the Lie derivative in the  
direction of the
gradient f. There are two types of  gradient  solitons on 
$\Bbb R^2$.  Namely,
the  flat soliton $(C_{\infty}=\infty, A = 1)$ and the 
cigar soliton
$(C_{\infty}<\infty, A = 0)$.  The flat soliton is the 
standard flat metric on
$\Bbb R^2$.  The cigar soliton is a metric, which can be 
expressed as $ds^2
=({du^2 + dv^2)/(1 + u^2 + v^2}),$ where $\{u, v\}$ are 
rectangular
coordinates on $\Bbb R^2.$ The cigar solitons are the  
so-called Barenblatt
solutions in the field of the porous medium equations.

The Ricci flow  and the soliton phenomenon gave a  new 
proof for the
uniformization theorem on compact surfaces and orbifolds 
without boundary
(\cite {C, Ha1,  W2, CW}).  Understanding  the 
solitons may provide insights toward studying the Ricci 
flow on 
higher-dimensional K\"ahler manifolds.  For an 
announcement of related work, see \cite
{Shi}.  As a step towards further studying the soliton 
phenomenon on 
higher-dimensional K\"ahler manifolds, Richard Hamilton 
raised  the following
question: 

\dfn{ Question}  On what manifolds do solutions to the 
Ricci flow 
asymptotically approach  nontrivial solitons\/{\rm ?} 
\enddfn
 
One of the simplest spaces is $\Bbb R^2$ with complete 
metric.  We say that 
the Ricci flow on $\Bbb R^2$  has {\it weak modified 
convergence\/} at time
infinity if  there exists a 1-parameter family of  
diffeomorphisms
$\{\phi_t\}_{t \in [0, \infty)}$ on $\Bbb R^2$ such that 
for any sequence of
times going to infinity there is a subsequence of times 
$\{ t_j \}_{j
=0}^{\infty}$ and the modified metric 
$ds^2(\phi_{t_j}(\,{\boldcdot}\,),t_j)$ 
converges uniformly on every compact set as $j \to 
\infty$.   We will  now
state our results concerning the Ricci flow.  
 
\thm{Main Theorem {\rm (Ricci flow \cite {W2})}\rm} Given  
a  complete $ (\Bbb
R^2, ds^2(0)) $ with $|R|\le C$ and  $ |Du| \le C$ at $t 
=0$, the Ricci
flow  has weak modified convergence  at time infinity to a 
limiting metric. In
the case when $ R>0$ at $t =0$, the limiting metric is a 
cigar soliton if
$C_{\infty}(ds^2(0)) <\infty$, or a flat metric if 
$A(ds^2(0)) >0$.
 \ethm

In  the process of proving the main theorem we have the 
following lemmas.
Let $R$ denote the scalar curvature and $ R_{-} = \max \{ 
-R, 0 \}.$

\thm{Lemma 1 {\rm (Ricci flow) (Long time existence)}\rm} 
Given  a  complete $ (
\Bbb R^2, ds^2) $ with $|R| \le C$  and $ |Du| \le C$  at 
$t=0$, under
the Ricci flow $(^*1)$,  the solution of the flow  exists 
for infinite time.
\ethm

\thm{Lemma 2 {\rm (Ricci flow)}\rm} Given  a  complete $ ( 
\Bbb R^2, ds^2) $ with
$|R| \le C$, $\int_{\Bbb R^2} R_{-} d\mu < \infty$, $ |Du| 
\le C$, and
$C_{\infty} > 0$  at $t=0$, under the Ricci flow $(^*1)$, 
we have the
following\/\,{\rm :}
 \roster
\item"{1.}" {\rm(Uniqueness)} The solution of the flow is 
unique. 
\item"{2.}" {\rm(Geometric Properties)} $C_{\infty}$, 
$A(g_{ij}) $,
and $\int_{\Bbb R^2} Rd\mu < \infty$ are constants
under the flow. 
\endroster
\ethm

\thm{ Lemma 3 {\rm (Ricci flow)}\rm} Given  a  complete $ 
( \Bbb R^2, ds^2) $ with
$0< R\le C$ and  $ |Du| \le C$ at $t=0$. Then, under the 
Ricci flow $(^*1)$,
$\lim_{ t \to \infty} e^{u(x,y,t)}$  converges uniformly 
on every compact
set and  $\lim_{ t \to \infty} e^{u(x,y,t)}$ is either 
identically zero or
positive everywhere.  If  $\lim_{ t \to \infty} 
e^{u(x,y,t)}>0$ then
$\lim_{ t \to \infty} e^{u(x,y,t)}(dx^2 + dy^2)$  induces 
a metric on
$\Bbb R^2$ with  curvature identically zero. \ethm

Nevertheless, it is possible  to choose a 1-parameter 
family of diffeomorphisms
to  get weak modified convergence.  To see why modifying 
the solution
by 
diffeomorphism is needed, we will illustrate the following 
example. 

\ex{Example 2.2} Given a cigar soliton (or one of the 
Barenblatt solutions) 
$ds^2(0)= ({dx^2 + dy^2 })/( {1 + x^2 + y^2})$ on $\Bbb 
R^2$, it is easy
to compute  that the solution of the Ricci flow  with 
initial data $ds^2(0)$ is 
$ds^2(t) = ({dx^2 + dy^2 })/ ({e^{4t} +x^2 + y^2}).$  Then 
$e^{u(x, y,t)} = 
1/( {e^{4t} +x^2 + y^2})$  goes to zero as  time 
approaches infinity; 
therefore, we cannot claim that $\lim_{t \to \infty} 
e^{u(x, y,t)}$ yields a
metric on $\Bbb R^2.$

Nevertheless, if we let diffeomorphism $\phi_t(A, B) = 
(e^{2t}A, e^{2t}B) =
(x,y)$, then 
$$ds^2(x,y,t)= ds^2(\phi_t(A, B), t)= \frac{dA^2 + dB^2 } 
{1 +A^2 + B^2}.$$
Let $\widehat{ds^2}(A,B, t) = ds^2(\phi_t(A, B), t)$ and 
$ e^{\widehat u} ={1 /}( {1+A^2 + B^2}).$  Then 
$e^{\widehat u}$ is 
stationary in time.
\endex
 
Note that the Ricci flow on other complete noncompact 
surfaces is also
discussed  in \cite {W2}.

Now we will list the corresponding relations between the 
function $v$ and the
geometric properties.  If $\{r, \theta \} $ are  polar 
coordinates on $\Bbb
R^2$, then $ ds^2 = v(dr^2 + r^2d\theta^2)$ and 
$$
\gather
(\Bbb R^2, ds^2) \text{ is complete }\Longrightarrow
 v>0 \text{  and  } \int_{r =0}^{\infty} 
v^{1/2}(\theta, r) \,dr = \infty \tag "(*4)"\\ 
\shoveright {\forall 0 \le \theta \le 2 \pi,}
\endgather$$
$$ |R| \le C\Longleftrightarrow
\left| \frac{\overline{\Delta}\ln v}{v}\right| \le C, \tag 
"(*5)" $$
$$ 0 < R \le C\Longleftrightarrow
 0 < - \frac{\overline{\Delta}\ln v}{v} \le C, \tag 
"$(^*5')$" $$
$$|Du| \le C\Longleftrightarrow |({ v_x^2 + v_y^2})/{
v^3}|\le C^2; \tag "(*6)"$$
$$\int R_{-} \,d\mu \le C\Longleftrightarrow 
\int\max \{ -\overline{\Delta} \ln v, 0 \}
\,dx \,dy \le C, \tag "(*7)"$$
$$ C_{\infty} = \lim_{r \to \infty}
\int rv^{1/2} d\theta > 0
\Longrightarrow v \notin L^1, \tag "(*8)"$$
$$A(ds^2) = 
\lim_{r \to \infty} \frac{\int rv^{1/2} d\theta}{\int 
v^{1/2} dr} > 0.
\tag "(*9)"
$$

The relations (*5), (*6), and (*8)  follow from 
$$ \gather
 \qquad R = - e^{-u}\overline{\Delta}u = -
\frac{\overline{\Delta}\ln v}{v} ,\tag "(**5)" \\
\qquad u = \ln v \quad \text{and} \quad |Du|^2 = \frac1v 
\frac{<v_x,
v_y>}{v}\,{\boldcdot}\, \frac{<v_x,v_y>}{v} , \tag "(**6)"\\
 \int v\,dx\,dy = \int e^u \,dx\,dy =\infty \.   \tag 
"(**8)"\\
 \endgather
$$

We say that  the PME $(m=0)$ on $\Bbb R^2$  has weak 
modified convergence 
at time infinity if  there exists a 1-parameter family of
reparametrizations 
$\{\phi_t\}_{t \in [0, \infty)}$ on $\Bbb R^2$ such that  
for any sequence of
times going to infinity there is a subsequence of times 
$\{ t_j \}_{j
=0}^{\infty}$ and the modified solution 
$v(\phi_{t_j}(\,{\boldcdot}\,),t_j)$ 
converges uniformly to a positive function on every 
compact set as $j \to
\infty$.   Then we have

\thm{ Main Theorem* {\rm [PME $(m=0)]$ }\rm} On $\Bbb 
R^2$, let the  positive
function $v$ satisfy $(^*4), (^*5)$, and  $(^*6)$   at 
$t=0$. Then, under the 
PME $(m=0)$, the solution $v(x,t)$ of PME $(m=0)$ with  
$v(x,0) =v(x)$
has weak modified convergence  at time infinity to a 
limiting positive function
$v_{\infty}$ satisfying $(^*4)$. In the case when 
$(*5^\prime)$ also holds at $
t=0$, $v_{\infty}$  is  one  of the Barenblatt solutions if
$C_{\infty}(v(\,{\boldcdot}\, , 0))< 
\infty$, or a constant if $A(v(\,{\boldcdot}\,,0)) >0$. 
\ethm

\thm{ Lemma *1 {\rm [PME $(m=0)]$  (Long time existence) 
}\rm} On $\Bbb R^2$, if
the  positive function $v$ satisfies $(^*4), (^*5)$, and  
$(^*6)$ at $t=0$, 
then, under the  PME $(m=0)$, the solution $v(x,t)$ of PME 
$(m=0)$ with 
$v(x,0) =v(x)$  exists for infinite time. \ethm

\thm{ Lemma *2 {\rm(PME $(m=0))$}\rm}\! \kern-.39pt On
\kern-.39pt $\Bbb R^2$\!,
\kern-.39pt if the \kern-.39pt positive \kern-.39pt 
function \kern-.39pt $v$
\kern-.39pt satisfies \kern-.39pt $(^*4)$,  
$(^*5)$,
$(^*6)$, $(^*7)$, and 
$(^*8)$  at $t=0$, 
then, under the  PME $(m=0)$,  we have the 
following\,\/{\rm :}
 \roster
\item"{1.}" {\rm(Uniqueness)} The solution of the flow is 
unique. 
\item"{2.}" {\rm(Geometric Properties)} $C_{\infty}$, 
$A(ds^2) $,
and $\int_{\Bbb R^2} R\,d\mu < \infty$ are constants
under the flow. 
\endroster
\ethm

\thm{ Lemma *3 {\rm (PME $(m=0))$ }\rm} On $\Bbb R^2$, if
the  positive function $v$ satisfies $(^*4),  
(^*5')$, and  $(^*6)$ at
$t=0$, then, under the  PME $(m=0)$,  $\lim_{ t \to \infty}
v(x,y,t)$   converges uniformly on every compact set and 
$\lim_{ t \to \infty} v(x,y,t)$ is either identically zero 
or positive
everywhere.  If  $\lim_{ t \to \infty} v(x,y,t)>0$ then 
$\lim_{ t \to
\infty} v(x,y,t)(dx^2 + dy^2)$  induces a metric on $\Bbb 
R^2$ with  curvature
identically zero, in particular, $\lim_{ t \to \infty} 
v(x,y,t)$ is a constant.
\ethm
 
Note that there is still a large class of Riemannian 
structures with
$C_{\infty} =\infty$ and  $A =0$, which our method fails 
to classify the limit.  

\demo{Sketch of the proof  } The evolution equation of  
$h= R +
|Du|^2$ provides the infinite time existence and 
uniform\ bounds for $|Du| $,
$| D^ku|$, $R$, and $|D^k R|$ for all $k \ge 1$ after a 
short time.  Finite
total curvature and  $C_{\infty} >0$ imply that the 
curvature  decays to zero
at distance infinity.  This yields that $C_{\infty}$, 
$A(ds^2)$, and $\int
R\,d\mu$ are preserved under the flow.  Furthermore, the 
solution of the flow is
unique and the injectivity radius $i(M)$ decays at most 
exponentially.

The positivity of the curvature of an initial metric 
provides pointwise
convergence of the function $e^{u}$ at time infinity.   
The uniform bounds on
$|D^m u|$ imply   $\lim_{ t \to \infty}e^{u}$ is a smooth 
function and is
either identically zero or positive everywhere.    In the 
case, when  $ \lim_{
t \to \infty}e^{u} >0$, the limiting solution is a flat 
metric.

We also may  choose a 1-parameter family of 
diffeomorphisms $\phi_t :  \Bbb R^2
\to \Bbb R^2$ such that there exists a sequence of times 
$\{ t_j \}_{j
=0}^{\infty}$ and $ \lim_{ t \to \infty} 
ds^2(\phi_{t_j}(\,{\boldcdot}\,),
t_j)$ converges uniformly on every compact set.  If $ R>0$ 
at $t=0$ and
$C_{\infty} (ds^2(0))<\infty,$ then some integral bound 
classifies $ \lim_{ t \to
\infty} ds^2(\phi_{t_j}(\,{\boldcdot}\,), t_j)$ as a cigar 
soliton with 
circumference no bigger than $C_{\infty} 
(ds^2(0))<\infty.$  If $ R>0$ at $t=0$ and
$A(ds^2(0)) >0$, then the  Harnack's inequality classifies 
$ \lim_{ t \to \infty}
ds^2(\phi_{t_j}(\,{\boldcdot}\,), t_j)$ as a  flat metric.
\enddemo

\enddocument